\documentclass[a4paper,11pt]{article}

\usepackage{amssymb}
\usepackage{amsmath}  
\usepackage{amsfonts} 

\begin{document}

\title{Asymptotics of the number of threshold functions on a
two-dimensional rectangular grid}

\author{Pentti Haukkanen and Jorma K. Merikoski\\
School of Information Sciences\\
FI-33014 University of Tampere, Finland\\\\
\texttt{pentti.haukkanen@uta.fi, jorma.merikoski@uta.fi}}

\maketitle

\begin{abstract}

Let $m,n\ge 2$, $m\le n$. It is well-known that the number of (two-dimensional) threshold functions on an $m\times n$ rectangular grid is
\begin{eqnarray*}
t(m,n)=\frac{6}{\pi^2}(mn)^2+O(m^2n\log{n})+O(mn^2\log{\log{n}})=
\\
\frac{6}{\pi^2}(mn)^2+O(mn^2\log{m}).
\end{eqnarray*}
We improve the error term by showing that
$$
t(m,n)=\frac{6}{\pi^2}(mn)^2+O(mn^2).
$$

\end{abstract}


\smallskip


\smallskip

{\it Running title:} Asymptotics of number of threshold
functions

\section{Introduction}

Consider a rectangular
grid $G=G(m,n)=\{0,\dots,m-1\}\times\{0,\dots,n-1\}$, where
$m,n\ge 2$. A function $\tau:G\to\{0,1\}$ is a
(two-dimensional) {\it threshold function} if there is a line
separating the sets
$\tau^{-1}(\{0\})$ and
$\tau^{-1}(\{1\})$. In other words, there are real numbers
$a$,
$b$ and $c$ such that $\tau(x,y)=0$ if and only if $ax+by+c\le
0$. Acketa and \v Zuni\'c~\cite[Theorem~3]{AZ} proved (with somewhat
different formulation) that the number of these functions is
\begin{eqnarray}
\label{tmnexpl}
t(m,n)=f(m,n)+2,
\end{eqnarray}
where 
\begin{eqnarray*}
\label{fmn}
f(m,n)=
\sum_{\substack{-m<i<m\\-n<j<n\\(i,j)=1}}(m-|i|)(n-|j|)
\end{eqnarray*}
and $(i,j)$ denotes the greatest common divisor of $i$
and~$j$. See also~\cite[p.~9--10]{Mu2}. In particular,
\begin{eqnarray}
\label{tn}
t(n)=t(n,n)=f(n)+2,
\end{eqnarray}
where
\begin{eqnarray*}
f(n)=
\sum_{\substack{-n<i,j<n\\(i,j)=1}}(n-|i|)(n-|j|).
\end{eqnarray*}

\smallskip

We are interested in asymptotic formulas. By Koplowitz et al.~\cite[Theorem~2]{KLB},
\begin{eqnarray}
\label{tnasympalex}
t(n)=\frac{6}{\pi^2}n^4+O(n^3\log{n}).
\end{eqnarray}
Let $m\le n$. By Acketa and \v Zuni\' c~\cite[p.~168]{AZ},
\begin{eqnarray}
\label{tmnasympaz}
t(m,n)=\frac{6}{\pi^2}(mn)^2+O(m^2n\log{n})+O(mn^2\log{\log{n}}).
\end{eqnarray}
By Alekseyev~\cite[Theorem~25]{Al},
\begin{eqnarray}
\label{tmnasympalex}
t(m,n)=\frac{6}{\pi^2}(mn)^2+O(mn^2\log{m}),
\end{eqnarray}
but see also \v Zuni\' c~\cite{Zu}. The error terms of (\ref{tmnasympaz}) and~(\ref{tmnasympalex}) are incomparable. For example, (\ref{tmnasympaz}) is better than~(\ref{tmnasympalex}) if $m=\sqrt{n}$, while (\ref{tmnasympalex}) is better if $m$ remains constant.

\smallskip

We will in Section~\ref{Improved formulas} improve the error terms of
these formulas. The improvements of~(\ref{tnasympalex})
follow directly from well-known results, but the improvement 
of (\ref{tmnasympaz}) and~(\ref{tmnasympalex}) must be proved; we will do so in
Section~\ref{Proof of}. Before it, we will in
Section~\ref{Basic facts and notations} introduce the basic
facts and notations needed in the proof. We will complete our
paper with conclusions and remarks in
Section~\ref{Conclusions and remarks}.

\section{Improved formulas}
\label{Improved formulas}

We can drop out the logarithm from~(\ref{tnasympalex}):
\begin{eqnarray}
\label{tnasympon3}
t(n)=\frac{6}{\pi^2}n^4+O(n^3).
\end{eqnarray}
We have also a sharper but more complicated formula
\begin{eqnarray}
\label{tnasymp}
t(n)=\frac{6}{\pi^2}n^4+O(n^3\exp(-A(\log{n})^\frac{3}{5}
(\log\log{n})^{-\frac{1}{5}}))
\end{eqnarray}
for some $A>0$. Assuming the Riemann hypothesis (RH), this
still improves into
\begin{eqnarray}
\label{tnasymprh}
t(n)=\frac{6}{\pi^2}n^4+O(n^{\frac{5}{2}+\varepsilon})
\end{eqnarray}
for all $\varepsilon>0$. To prove (\ref{tnasymp}) and
(\ref{tnasymprh}), we simply note that $t(n)$ and $f(n)$ have by~(\ref{tn}) the same asymptotic behavior and refer to~\cite[Lemma~6]{EMHM}.

\smallskip

We can also drop
out the logarithms from (\ref{tmnasympaz}) and~(\ref{tmnasympalex}).
Generalizing~(\ref{tnasympon3}), we will in
Section~\ref{Proof of} prove that
\begin{eqnarray}
\label{tmnasymp}
t(m,n)=\frac{6}{\pi^2}(mn)^2+O(mn^2)
\end{eqnarray}
if $m\le n$. We conjecture~\cite{Mu3} that the error term
can be improved by substituting $n\mapsto\sqrt{mn}$ in
(\ref{tnasymp}) and~(\ref{tnasymprh}). That is,
\begin{eqnarray}
\label{tmnasympconj}
t(m,n)=\frac{6}{\pi^2}(mn)^2+O((mn)^\frac{3}{2}\exp(-A(\log{mn})^\frac{3}{5}
(\log\log{mn})^{-\frac{1}{5}}))
\end{eqnarray}
for some $A>0$ and, under~RH,
\begin{eqnarray}
\label{tmnasymprhconj}
t(m,n)=\frac{6}{\pi^2}(mn)^2+O((mn)^{\frac{5}{4}+\varepsilon})
\end{eqnarray}
for all $\varepsilon>0$.

\section{Basic facts and notations}
\label{Basic facts and notations}

The ''delta function''~$\delta$,
defined by
\begin{eqnarray*}
\delta(k)=\left\{\begin{array}{ccc}
1&\mathrm{if}&k=
1,\\0&\mathrm{if}&k\ge 2,
\end{array}\right\}
\end{eqnarray*}
and the M\"obius function~$\mu$, defined by
\begin{eqnarray*}
\mu(k)=\left\{\begin{array}{cc}1&\mathrm{if}\,\,k=1,\qquad
\qquad\qquad\qquad\qquad\qquad\\
0&\mathrm{if}\,\,k\,\,\mathrm{has\,one\,or\,more\,repeated\,
prime\,factors},\\(-1)^r&\mathrm{if}
\,\,k\,\,\mathrm{is\,a\,product\,of}\,\,r\,\,
\mathrm{distinct\,primes},
\end{array}\right\}
\end{eqnarray*}
satisfy
\begin{eqnarray}
\label{deltamu}
\delta(k)=\sum_{d\,|\,k}\mu(d).
\end{eqnarray}

We give some asymptotic formulas involving~$\mu(k)$. Let
$0\le c_1,\dots,c_m\le 1$. Since
\begin{eqnarray}
\label{cmu}
c_k|\mu(k)|\le 1\quad(k=1,\dots,m),
\end{eqnarray}
we have
\begin{eqnarray}
\label{mu1}
\sum_{k=1}^mc_k\mu(k)=O(m).
\end{eqnarray}
We also have
\begin{eqnarray}
\label{mu2}
\sum_{k=1}^m\frac{\mu(k)}{k}=O(1),
\end{eqnarray}
which in fact can be sharpened by replacing $O$
with~$o$~\cite[p.~194]{MSC}. Furthermore, (\ref{cmu}) implies easily 
that
\begin{eqnarray}
\label{mu3}
\sum_{k=1}^mc_k\frac{\mu(k)}{k}=O(\log{m}).
\end{eqnarray}
Finally, we give
\begin{eqnarray}
\label{mu4}
\sum_{k=1}^m\frac{\mu(k)}{k^2}=\frac{6}{\pi^2}+
O\big(\frac{1}{m}\big),
\end{eqnarray}
which in fact can be sharpended by substituting on the
right-hand side $m\mapsto m(\log{m})^a$ for
any~$a>0$~\cite[p.~194]{MSC}.

\smallskip

Let us define
$$
\{x\}=x-\lfloor
x\rfloor,\qquad\alpha(m,n)=
\sum_{d=1}^m\frac{\mu(d)}{d}\big\{\frac{n}{d}\big\},\qquad
\alpha(m)=\alpha(m,m).
$$
Because $0\le\{x\}<1$, we see by~(\ref{mu3}) that
\begin{eqnarray}
\label{gamma}
\alpha(m,n)=O(\log{m}).
\end{eqnarray}

\section{Proof of (\ref{tmnasymp})}
\label{Proof of}

By~(\ref{tmnexpl}), an equivalent task is to prove that
\begin{eqnarray}
\label{fmnasymp}
f(m,n)=\frac{6}{\pi^2}(mn)^2+O(mn^2)
\end{eqnarray}
if $m\le n$. We divide the proof into six
parts.

\medskip
\noindent
{\it 1. Evaluating} $f(m+1,n+1)$. We have
\begin{eqnarray*}
 f(m+1,n+1)=\sum_{\substack{-m\le
i\le m\\-n\le j\le n\\{(i,j)}=1}}(m+1-|i|)
(n+1-|j|)=
\\
4\sum_{\substack{1\le i \le m\\1\le j\le n\\(i,j)=1}}(m+1-i)
(n+1-j)+2(m+1)n+2m(n+1)=
\end{eqnarray*}
\begin{eqnarray}
\nonumber
4\sum_{\substack{1\le
i\le m\\1\le j\le n\\(i,j)=1}}\left[(m+1)(n+1)-(n+1)i-(m+1)j
+ij\right]+4mn+2m+2n
\\\nonumber
=4\Big[(m+1)(n+1)\sum_{\substack{1\le i\le m\\1\le j\le
n\\(i,j)=1}} 1-(n+1)
\sum_{\substack{1\le i\le m\\i\le j\le n\\(i,j)=1}} i-(m+1)
\sum_{\substack{1\le i\le m\\1\le j\le n\\(i,j)=1}}
j+
\\\nonumber
\sum_{\substack{1\le i\le m\\1\le
j\le n\\(i,j)=1}} ij\Big]+O(mn)=
4\big[(m+1)(n+1)s_1(m,n)-(n+1)s_2(m,n)-
\\
(m+1)s_3(m,n)+
\label{ffmn}
s_4(m,n)\big]+O(mn), 
\end{eqnarray}
where
\begin{eqnarray*}
s_1(m,n)=\sum_{\substack{1\le i\le m\\1\le j\le n\\(i,j)=1}}
1,\qquad s_2(m,n)=\sum_{\substack{1\le i\le
m\\1\le j\le n\\(i,j)=1}} i,
\\s_3(m,n)=\sum_{\substack{1\le i\le m\\1\le j\le n\\(i,j)=1}}
j,\qquad s_4(m,n)=\sum_{\substack{1\le i\le m\\1\le j\le
n\\(i,j)=1}} ij.
\end{eqnarray*}

\medskip
\noindent
{\it 2. Evaluating} $s_1(m,n)$. Remember that $m\le n$.
Denote
$i=ad$, $j=bd$ and apply~(\ref{deltamu}). Then,
by~(\ref{mu1}), (\ref{mu2}), (\ref{mu4}) and~(\ref{gamma}),
\begin{eqnarray*}
\label{s1}
s_1(m,n)=\sum_{i=1}^m\sum_{j=1}^n\delta((i,j))=
\sum_{i=1}^m\sum_{j=1}^n\sum_{d|(i,j)}\mu(d)=
\sum_{d=1}^m\mu(d)\sum_{a=1}^{\lfloor\frac{m}{d}\rfloor}
\sum_{b=1}^{\lfloor\frac{n}{d}\rfloor}1=
\\
\sum_{d=1}^m\mu(d)\lfloor\frac{m}{d}\rfloor
\lfloor\frac{n}{d}\rfloor=\sum_{d=1}^m\mu(d)\big(\frac{m}{d}-
\big\{\frac{m}{d}\big\}\big)\big(\frac{n}{d}-
\big\{\frac{n}{d}\big\}\big)=\qquad
\end{eqnarray*}
\begin{eqnarray}
\nonumber
mn\sum_{d=1}^m\frac{\mu(d)}{d^2}-
m\sum_{d=1}^m\frac{\mu(d)}{d}\big\{\frac{n}{d}\big\}
-n\sum_{d=1}^m\frac{\mu(d)}{d}\big\{\frac{m}{d}\big\}
+\sum_{d=1}^m\mu(d)
\big\{\frac{m}{d}\big\}\big\{\frac{n}{d}\big\}=
\nonumber\\
mn\big(\frac{6}{\pi^2}+O\big(\frac{1}{m})\big)-m\alpha(m,n)
-n\alpha(m)+O(m)=\qquad\qquad\qquad\qquad
\nonumber\\
\frac{6}{\pi^2}mn-m\alpha(m,n)
-n\alpha(m)+O(n).\qquad
\end{eqnarray}

\medskip
\noindent
{\it 3. Evaluating} $s_2(m,n)$. Similarly,
\begin{eqnarray*}
2s_2(m,n)=2\sum_{i=1}^m\sum_{j=1}^n
i\delta((i,j))=2\sum_{i=1}^m\sum_{j=1}^n
i\sum_{d|(i,j)}\mu(d)=\qquad\qquad
\\
2\sum_{d=1}^m\mu(d)d\sum_{a=1}^{\lfloor\frac{m}{d}\rfloor}
\sum_{b=1}^{\lfloor\frac{n}{d}\rfloor}a=
2\sum_{d=1}^m\mu(d)d\,\frac{1}{2}\lfloor\frac{m}{d}\rfloor
\big(\lfloor\frac{m}{d}\rfloor+1\big)
\lfloor\frac{n}{d}\rfloor=\qquad
\\
\sum_{d=1}^m\mu(d)d\big(\frac{m}{d}-
\big\{\frac{m}{d}\big\}\big)\big(\frac{m}{d}-
\big\{\frac{m}{d}\big\}+1\big)\big(\frac{n}{d}-
\big\{\frac{n}{d}\big\}\big)=\qquad
\\
m^2n\sum_{d=1}^m\frac{\mu(d)}{d^2}+
mn\sum_{d=1}^m\frac{\mu(d)}{d}
-m^2\sum_{d=1}^m\frac{\mu(d)}{d}\big\{\frac{n}{d}\big\}
-2mn\sum_{d=1}^m\frac{\mu(d)}{d}\big\{\frac{m}{d}\big\}
\\
-m\sum_{d=1}^m\mu(d)\big\{\frac{n}{d}\big\}
-n\sum_{d=1}^m\mu(d)\big\{\frac{m}{d}\big\}
+n\sum_{d=1}^m\mu(d)\big\{\frac{m}{d}\big\}^2+\qquad\qquad
\\
2m\sum_{d=1}^m\mu(d)\big\{\frac{m}{d}\big\}
\big\{\frac{n}{d}\big\}+
\sum_{d=1}^m\mu(d)d\big\{\frac{m}{d}\big\}
\big\{\frac{n}{d}\big\}-
\sum_{d=1}^m\mu(d)d\big\{\frac{m}{d}\big\}^2\big\{\frac{n}{d}\big\}=
\\
m^2n\big(\frac{6}{\pi^2}+O(\frac{1}{m})\big)+mnO(1)-
m^2\alpha(m,n)-2mn\alpha(m)+\qquad
\\
mO(m)+
nO(m)+nO(m)+mO(m)+O(m^2)+O(m^2)=\qquad
\\
\frac{6}{\pi^2}m^2n-m^2\alpha(m,n)-2mn\alpha(m)+O(mn),\qquad
\end{eqnarray*}
and so
\begin{eqnarray}
\label{s2}
s_2(m,n)=\frac{3}{\pi^2}m^2n-\frac{1}{2}m^2\alpha(m,n)-
mn\alpha(m)+
O(mn).
\end{eqnarray}

\medskip
\noindent
{\it 4. Evaluating} $s_3(m,n)$. A simple modification of
Part~3 gives
\begin{eqnarray}
\label{s3}
s_3(m,n)=\frac{3}{\pi^2}mn^2-mn\alpha(m,n)-
\frac{1}{2}n^2\alpha(m)+O(mn).
\end{eqnarray}

\medskip
\noindent
{\it 5. Evaluating} $s_4(m,n)$. Now we need quite heavy
calculation. We obtain
\begin{eqnarray*}
4s_4(m,n)=4\sum_{i=1}^m\sum_{j=1}^n
ij\delta((i,j))=4\sum_{i=1}^m\sum_{j=1}^n
ij\sum_{d|(i,j)}\mu(d)=\qquad\qquad
\\
4\sum_{d=1}^m\mu(d)d^2\sum_{a=1}^{\lfloor\frac{m}{d}\rfloor}
\sum_{b=1}^{\lfloor\frac{n}{d}\rfloor}ab=4\sum_{d=1}^m\mu(d)d^2\,\frac{1}{2}\lfloor\frac{m}{d}\rfloor
\big(\lfloor\frac{m}{d}\rfloor+1\big)
\frac{1}{2}\lfloor\frac{n}{d}\rfloor\big(\lfloor\frac{n}{d}\rfloor+1\big)
=
\\
\sum_{d=1}^m\mu(d)d^2\big(\frac{m}{d}-
\big\{\frac{m}{d}\big\}\big)\big(\frac{m}{d}-
\big\{\frac{m}{d}\big\}+1\big)\big(\frac{n}{d}-
\big\{\frac{n}{d}\big\}\big)\big(\frac{n}{d}-
\big\{\frac{n}{d}\big\}+1\big)=
\\
m^2n^2\sum_{d=1}^m\frac{\mu(d)}{d^2}+
m^2n\sum_{d=1}^m\frac{\mu(d)}{d}+
mn^2\sum_{d=1}^m\frac{\mu(d)}{d}-
2m^2n\sum_{d=1}^m\frac{\mu(d)}{d}\big\{\frac{n}{d}\big\}-
\\
2mn^2\sum_{d=1}^m\frac{\mu(d)}{d}\big\{\frac{m}{d}\big\}-
m^2\sum_{d=1}^m\mu(d)\big\{\frac{n}{d}\big\}-
n^2\sum_{d=1}^m\mu(d)\big\{\frac{m}{d}\big\}+\qquad\qquad
\\
mn\sum_{d=1}^m\mu(d)+
m^2\sum_{d=1}^m\mu(d)\big\{\frac{n}{d}\big\}^2+
n^2\sum_{d=1}^m\mu(d)\big\{\frac{m}{d}\big\}^2-
2mn\sum_{d=1}^m\mu(d)\big\{\frac{m}{d}\big\}-
\\
2mn\sum_{d=1}^m\mu(d)\big\{\frac{n}{d}\big\}+
4mn\sum_{d=1}^m\mu(d)\big\{\frac{m}{d}\big\}
\big\{\frac{n}{d}\big\}-
m\sum_{d=1}^m\mu(d)d\big\{\frac{n}{d}\big\}-\qquad
\\
n\sum_{d=1}^m\mu(d)d\big\{\frac{m}{d}\big\}+
2m\sum_{d=1}^m\mu(d)d\big\{\frac{m}{d}\big\}
\big\{\frac{n}{d}\big\}
+2n\sum_{d=1}^m\mu(d)d\big\{\frac{m}{d}\big\}
\big\{\frac{n}{d}\big\}+
\\
m\sum_{d=1}^m\mu(d)d\big\{\frac{n}{d}\big\}^2+
n\sum_{d=1}^m\mu(d)d\big\{\frac{m}{d}\big\}^2-
2m\sum_{d=1}^m\mu(d)d\big\{\frac{m}{d}\big\}
\big\{\frac{n}{d}\big\}^2-\qquad
\\
2n\sum_{d=1}^m\mu(d)d\big\{\frac{m}{d}\big\}^2
\big\{\frac{n}{d}\big\}+
\sum_{d=1}^m\mu(d)d^2\big\{\frac{m}{d}\big\}
\big\{\frac{n}{d}\big\}-\sum_{d=1}^m\mu(d)d^2
\big\{\frac{m}{d}\big\}^2
\big\{\frac{n}{d}\big\}-
\\
\sum_{d=1}^m\mu(d)d^2\big\{\frac{m}{d}\big\}
\big\{\frac{n}{d}\big\}^2+
\sum_{d=1}^m\mu(d)d^2\big\{\frac{m}{d}\big\}^2
\big\{\frac{n}{d}\big\}^2=
\\
m^2n^2\big(\frac{6}{\pi^2}+O(\frac{1}{m})\big)+m^2nO(1)+
mn^2O(1)-2m^2n\alpha(m,n)-2mn^2\alpha(m)+
\\
m^2O(m)+n^2O(m)+mnO(m)+m^2O(m)+n^2O(m)+mnO(m)+\qquad
\\
mnO(m)+mnO(m)+mO(m^2)+nO(m^2)+mO(m^2)+nO(m^2)+\qquad
\\
mO(m^2)+nO(m^2)+mO(m^2)+nO(m^2)+O(m^3)+O(m^3)+O(m^3)+
\\
O(m^3)=\frac{6}{\pi^2}m^2n^2-2m^2n\alpha(m,n)-2mn^2\alpha(m)+O(mn^2).
\end{eqnarray*}
Hence
\begin{eqnarray}
\label{s4}
s_4(m,n)=\frac{3}{2\pi^2}(mn)^2-\frac{1}{2}m^2n\alpha(m,n)-
\frac{1}{2}mn^2\alpha(m)+O(mn^2).
\end{eqnarray}

\medskip
\noindent
{\it 6. Final computation}. By (\ref{ffmn}), (\ref{s1}),
(\ref{s2}), (\ref{s3}), (\ref{s4}) and~(\ref{gamma}),
\begin{eqnarray*}
a=(m+1)(n+1)s_1(m,n)=\qquad\qquad\qquad\qquad\qquad\qquad\qquad\qquad\qquad
\\
(m+1)(n+1)\big(\frac{6}{\pi^2}mn-m\alpha(m,n)
-n\alpha(m)+O(n)\big)=
\\
\frac{6}{\pi^2}(mn)^2+\frac{6}{\pi^2}m^2n+\frac{6}{\pi^2}mn^2
+
\frac{6}{\pi^2}mn-m^2n\alpha(m,n)-mn^2\alpha(m)-
\\
m^2\alpha(m,n)-
n^2\alpha(m)-mn\alpha(m,n)-mn\alpha(m)-m\alpha(m,n)-\qquad
\\
n\alpha(m)+O(mn^2)=
\\
\frac{6}{\pi^2}(mn)^2-m^2n\alpha(m,n)-mn^2\alpha(m)+O(mn^2),
\end{eqnarray*}
\begin{eqnarray*}
b=(n+1)s_2(m,n)=(n+1)\big(\frac{3}{\pi^2}m^2n-
\frac{1}{2}m^2\alpha(m,n)-\qquad\qquad\qquad
\\
mn\alpha(m)+
O(mn)\big)=
\\
\frac{3}{\pi^2}(mn)^2+\frac{3}{\pi^2}m^2n-
\frac{1}{2}m^2n\alpha(m,n)
-mn^2\alpha(m)+O(mn^2),
\end{eqnarray*}
\begin{eqnarray*}
c=(m+1)s_3(m,n)=(m+1)\big(\frac{3}{\pi^2}mn^2-mn\alpha(m,n)-
\qquad\qquad\qquad
\\
\frac{1}{2}n^2\alpha(m)+O(mn)\big)=
\\
\frac{3}{\pi^2}(mn)^2+\frac{3}{\pi^2}mn^2-m^2n\alpha(m,n)
-\frac{1}{2}mn^2\alpha(m)+O(m^2n),
\end{eqnarray*}
\begin{eqnarray*}
d=s_4(m,n)=\frac{3}{2\pi^2}(mn)^2-\frac{1}{2}m^2n\alpha(m,n)-
\frac{1}{2}mn^2\alpha(m)+O(mn^2),\quad
\end{eqnarray*}
\begin{eqnarray*}
f(m+1,n+1)=4(a-b-c+d)+O(mn)=\frac{6}{\pi^2}(mn)^2+O(mn^2),
\end{eqnarray*}
which implies also (\ref{fmnasymp}).

\section{Conclusions and remarks}
\label{Conclusions and remarks}

We improved (\ref{tnasympalex}) into~(\ref{tnasympon3}) and
further into~(\ref{tnasymp}). Assuming~RH, we still
improved (\ref{tnasymp}) into~(\ref{tnasymprh}). An
interesting converse problem arises (cf. \cite[p.~168]{EMHM}):
Does (\ref{tnasymprh}) imply RH?

\smallskip

We also improved (\ref{tmnasympaz}) and~(\ref{tmnasympalex}) into~(\ref{tmnasymp}).
Because of~(\ref{tnasymp}) and~(\ref{tnasymprh}), we
conjectured that (\ref{tmnasymp}) can be further improved
into~(\ref{tmnasympconj}) and, assuming~RH, still
into~(\ref{tmnasymprhconj}).

\smallskip

More generally, given $q\ge 1$, define
$$
f_q(m,n)=
\sum_{\substack{-m<i<m\\-n<j<n\\{(i,j)}=q}}(m-|i|)(n-|j|).
$$
The proof of~(\ref{fmnasymp}) can be extended
to show that
$$
f_q(m,n)=\frac{6}{\pi^2q^2}(mn)^2+O(mn^2)
$$
if $m\le n$. Also the formulas (\ref{tnasympon3}),
(\ref{tnasymp}) and (\ref{tnasymprh}) for~$f(n)$ generalize.
We have~\cite{HM2}
$$
f_q(n)=\frac{6}{\pi^2q^2}n^4+r(n),
$$
where $r(n)$ has the $O$-estimates given in the original
formulas.

\smallskip

Besides the number of threshold functions, there are several other quantities whose asymptotics can be studied in a similar way. For example, asymptotic formulas for the number of gridlines~\cite{EMHM} and $q$-gridlines~\cite{HM2} are well-known in~$G(n)$. A simple modification of our above procedure gives such formulas in~$G(m,n)$.

\smallskip

In practical computation of all these quantities, recursive formulas \cite{Mu1, Mu2,EMHM,HM1} are useful. They have also been applied in computer experiments to find asymptotic formulas.

\section*{Acknowledgment}

We thank the referees. One's suggestions improved the presentation. The other alerted to us references we were not aware of.

\end{document}